\input ./style/arxiv-general.cfg
\documentclass[MSNbibl,nameyear,secthm,seceqn,dvips]{arxstspdf}
\makeatletter
   \@ifpackageloaded{graphicx}{}{\usepackage{graphicx}}
\makeatother
\usepackage[mathscr]{euscript}
\usepackage{flushend}
\usepackage{stfloats}
\usepackage{url,breakurl}

\volume{30}
\issue{4}
\pubyear{2015}
\firstpage{518}
\lastpage{532}
\doi{10.1214/15-STS526}
\docsubty{FLA}

\makeatletter
\newtheorem{theorem}[thm]{Theorem}
\newcommand{\bx}{\mathbf{x}}
\newcommand{\bmu}{\bolds{\mu}}
\newcommand{\bgamma}{\bolds{\gamma}}
\newcommand{\bnu}{\bolds{\nu}}
\newcommand{\bX}{\mathbf{X}}
\newcommand{\D}{\mathsf{D}}
\renewcommand{\citep}[1]{\citeauthor{#1}, \citeyear{#1}}
\makeatother

\begin{document}
\begin{frontmatter}

\title{A Population Background for Nonparametric Density-Based
Clustering}
\runtitle{Nonparametric Density-Based Clustering}

\begin{aug}
\author[A]{\fnms{Jos\'e E.}~\snm{Chac\'{o}n}\corref{}\ead[label=e1]{jechacon@unex.es}}
\runauthor{J. E. Chac\'on}

\affiliation{Universidad de Extremadura}

\address[A]{Jos\'e E. Chac\'{o}n is Profesor Titular (Associate Professor),
Departamento de Matem\'{a}ticas,
Universidad de Extremadura,
06006 Badajoz,
Spain \printead{e1}.}
\end{aug}

\begin{abstract}
Despite its popularity, it is widely recognized that the investigation
of some theoretical aspects of clustering has been relatively sparse.
One of the main reasons for this lack of theoretical results is surely
the fact that, whereas for other statistical problems the theoretical
population goal is clearly defined (as in regression or
classification), for some of the clustering methodologies it is
difficult to specify the population goal to which the data-based
clustering algorithms should try to get close. This paper aims to
provide some insight into the theoretical foundations of clustering by
focusing on two main objectives: to provide an explicit formulation for
the ideal population goal of the modal clustering methodology, which
understands clusters as regions of high density; and to present two new
loss functions, applicable in fact to any clustering methodology, to
evaluate the performance of a data-based clustering algorithm with
respect to the ideal population goal. In particular, it is shown that
only mild conditions on a sequence of density estimators are needed to
ensure that the sequence of modal clusterings that they induce is consistent.
\end{abstract}

\begin{keyword}
\kwd{Clustering consistency}
\kwd{distance in measure}
\kwd{Hausdorff distance}
\kwd{modal clustering}
\kwd{Morse theory}\vspace*{-6pt}
\end{keyword}
\end{frontmatter}

\section{Introduction}\label{sec1}

Clustering is one of the branches of Statistics with more research
activity in recent years. As
noted by \citet{M07}, ``clustering is a young domain of research, where
rigorous methodology is
still striving to emerge.'' Indeed, some authors have recently expressed
their concerns about the
lack of theoretical or formal developments for clustering, as, for
instance, \citet{vLBD05}, \citet{BLP06}, \citet{ABD08}, \citet{ZBD09}. This paper
aims to contribute to this regularization (or, say, rigorousization).

Stated in its most simple form, cluster analysis consists in
``partitioning a data set into groups
so that the points in one group are similar to each other and are as
different as possible from the
points in other groups'' (\citep{HMS01}, page 293). Posed as such, the
problem does not even seem to
have a statistical meaning. In fact, in concordance with \citet{LRL07},
it is possible to roughly
classify clustering methods into three categories, depending on the
amount of statistical
information that they involve. These categories are very basically
depicted in the following three paragraphs.

Some clustering techniques are solely based on the distances between
the observations. Close observations are joined together to form a
group, and extending the notion of inter-point distance to distance
between groups, the resulting groups are gradually merged until all the
initial observations are contained into a single group. This
represents, of course, the notion of agglomerative \textit{hierarchical
clustering} (\citep{I08}, Section~12.3). The graphical outcome depicting
the successive agglomeration of data points up to a single group is the
well-known dendrogram, and depending on the notion of inter-group
distance used along the merging process, the most common procedures of
this type are known as single linkage, complete linkage or average
linkage (see also \citep{HTF09},  page~523).

A first statistical flavor is noticed when dealing with those
clustering methodologies that
represent each cluster by a central point, such as the mean, the median
or, more generally, a
trimmed mean. This class of techniques is usually referred to as \textit{partitioning methods}, and
surely the most popular of its representatives is $K$-means (\citep{MQ67}). For a prespecified
number $K$ of groups, these algorithms seek for $K$ centers with the
goal of optimizing a certain
score function representing the quality of the clustering
(\citep{Eal11}, Chapter~5).

When a more extended set of features of the data-generating probability
distribution is used to determine the clustering procedure, it is usual
to refer to these techniques as distribution-based clustering or, for
the common case of continuous distributions, as \textit{density-based
clustering}. This approach is strongly supported by some authors, like
\citet{CM13}, who explicitly state that ``density needs to be
incorporated in the clustering procedures.''

As with all the statistical procedures, there exist parametric and
nonparametric methodologies for density-based clustering.
Surely the gold standard of parametric density-based clustering is
achieved through mixture modeling, as
clearly described in \citet{FR02}. It is assumed that the distribution
generating the data is a
mixture of simple parametric distributions, for example, multivariate
normal distributions, and each
component of the mixture is associated to a different population
cluster. Maximum likelihood is
used to fit a mixture model and then each data point is assigned to the
most likely component using
the Bayes rule.

The nonparametric methodology is based on identifying clusters as
regions of high density separated
from each other by regions of lower density (\citep{W69}, \citep{H75}). Thus, a
cluster is seen as a zone of
concentration of probability mass. In this sense, population clusters
are naturally associated with
the modes (i.e., local maxima) of the probability density function, and
this nonparametric approach
is denominated mode-based clustering or \textit{modal clustering} (\citep{LRL07}). Precisely, each
cluster is usually understood as the ``domain of attraction'' of a mode
(\citep{S03}).

\begin{figure}[b]

\includegraphics{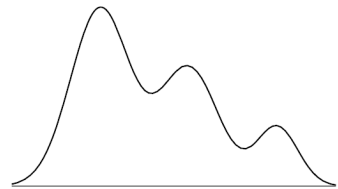}

\caption{Univariate trimodal density for which it is not possible to
capture its whole cluster structure using a level set analysis based on
a single level.}
\label{fig1}
\end{figure}

The concept of domain of attraction is not that simple to specify, and
providing a precise definition for that is one of the main goals of
this paper.
The first attempt to make the goal of modal clustering precise was
introduced through the notion of
level sets (\citep{H75}). If the distribution of the data has a
density~$f$, given $\lambda\geq0$, the
$\lambda$-level set of $f$ is defined as $L(\lambda)=\{\bx\colon f(\bx)\geq\lambda\}$. Then,
population $\lambda$-clusters are defined as the connected components
of $L(\lambda)$, a definition
that clearly captures the notion of groups having a high density. An
extensive account of the
usefulness of level sets in applications is given in \citet{MP09}.

One of the advantages of clustering based on level sets is that the
population target is clearly identified (the connected components of
the $\lambda$-level set). However, the main drawback of this approach
is perhaps the fact that the notion of population cluster
depends on the level $\lambda$, as recognized by \citet{S03}.
Nevertheless, other authors, like
\citet{CFF01} or \citet{CPP11}, affirm that the choice of $\lambda$ is
only a matter of resolution
level of the analysis.

Still, it is easy to think of many examples in which it is impossible
to observe the whole
cluster structure on the basis of a single level $\lambda$. Essentially
as in \citet{RSNW12}, page 906, Figure~\ref{fig1} shows a simple univariate example of
this phenomenon: three
different modal groups are visually identifiable, yet none of the level
sets of the density has three
connected components. To amend this, the usual recommendation is to
analyze the cluster structure
for several values of the level $\lambda$. Graphical tools oriented to
this goal are the cluster tree (\citep{S03}) or the mode
function (\citep{AT07}, \citep{MA14}). Both
graphics are useful to show how the clusters emerge as a function of
$\lambda$. See Section~\ref{Morse} for a more detailed explanation.

Finally, the idea of examining the evolution of the cluster structure
as the density level varies is closely related with the topic of
persistent homology, a tool from Computational Topology that, since its
relatively recent introduction, has attracted a great deal of interest
for its applications in Topological Data Analysis; see \citet{EH08},
\citet{C09} or \citet{CGOS13}. This tool allows to quantify which
topological aspects of an object are most persistent as the resolution
level evolves, thus leading to the identification of the most important
features of the object. In the context of data-based clustering based
on level sets, it can be very useful to distinguish which of the
discovered clusters are real and which of them are spurious (\citep{Fal14}).

The rest of this paper is structured as follows: in Section~\ref{sec2}
we introduce the concept of whole-space clustering as the type of
object of interest in cluster analysis, and we point out the difference
with the more usual notion of a clustering of the data. Later, it is
explained that the population whole-space clustering depends on the
adopted definition of cluster for each of the clustering methodologies.
Section~\ref{Morse} expands on the first main contribution of the paper
by providing a precise definition of the population goal of modal
clustering, making use of Morse theory, leading to an equivalent yet
simpler formulation (in a sense) as with the cluster tree. Once a
population background for clustering has been set up, Section~\ref{seccons} contains the second main contribution of the paper, a
proposal of two new loss functions to measure the similarity of two
whole-space clusterings. These distance functions are not limited to
modal clustering nor even to density-based clustering, they are
applicable to any clustering methodology having a clearly identified
population goal. As such, they can be used to define a notion of \textit{clustering consistency}, and for the particular case of modal
clustering it is shown that mild conditions are needed so that the
data-based clustering constructed from a sequence of density estimators
is consistent in this sense.

\section{Population Clusterings}\label{sec2}

Many different notions of cluster are possible, but no matter which one
is used, it is necessary to have a clear idea of the type of object
that clustering methods pursue from a population point of view. That
object will be called a clustering.

Since the empirical formulation of the clustering task comprises
partitioning a data set into groups, it suggests that its population
analogue should involve a partition of the whole space or, at least, of
the support of the distribution. Hence, a clustering of a probability
distribution $P$ on $\mathbb{R}^d$, or a \textit{whole-space $P$-clustering}, should be understood as an essential \mbox{partition}
of $\mathbb{R}^d$ into mutually disjoint measurable components, each
with positive
probability content (\citep{BLP06}). More specifically, a whole-space
$P$-clustering (or, simply, a clustering) is defined as a class
of measurable sets $\mathscr{C}=\{C_1,\ldots,C_r\}$ such that:
\begin{longlist}[3.]
\item[1.] $P(C_i)>0$ for all $i=1,\dots,r$,
\item[2.] $P(C_i\cap
C_j)=0$ for $i\neq j$, and
\item[3.] $P(C_1\cup\cdots\cup C_r)=1$.
\end{longlist}
The components $C_1,\ldots,C_r$ of such a partition are
called clusters. Thus, two clusterings $\mathscr{C}$ and $\mathscr{D}$ are
identified to be the same if they have the same number of clusters and,
up to a permutation of the
cluster labels, every cluster in $\mathscr{C}$ and its most similar
match in $\mathscr{D}$ differ in
a null-probability set (more details on this are elaborated in Section~\ref{seccons}).

At this point it is worth distinguishing between two different,
although closely related, concepts.
When the probability distribution $P$ is unknown, and a sample drawn
from $P$ is given, any
procedure to obtain a data-based (essential) partition
$\widehat{\mathscr{C}}=\{\widehat{C}_1,\ldots,\widehat{C}_r\}$ will be called a
\textit{data-based
clustering}. This simply means that $\int_{\widehat{C}_i}\,dP>0$ for all
$i=1,\ldots,r$, $\int_{\widehat{C}_i\cap
\widehat{C}_j}\,dP=0$ for $i\neq j$ and $\int_{\widehat{C}_1\cup\cdots\cup
\widehat{C}_r}\,dP=1$. 
However, when data are available most clustering procedures focus on
partitioning the data set, and, indeed, many of them do not even induce
a clustering of the probability distribution. This will be referred to
henceforth as a \textit{clustering of the data}. Notice that, clearly, any
data-based clustering $\widehat{\mathscr{C}}=\{\widehat{C}_1,\ldots
,\widehat{C}_r\}$ immediately results in a clustering of the data, by
assigning the same group to data points belonging to the same component
in $\widehat{\mathscr{C}}$.

\subsection{The Ideal Population Clustering}

The definition of (whole-space) clustering represents the type of
population object that clustering methods should
try to get close in general, but it is the particular employed notion
of cluster that makes the theoretical goal of
clustering methodologies change, focusing on different concepts of \textit{ideal population clustering}.

For some clustering techniques, this ideal population clustering is
well established. For instance,
it is well known that the population clustering induced by the optimal
set of $K$-means is a Voronoi
tessellation. To be precise, let $\bolds{\mu}_1^*,\ldots,\bolds{\mu}_K^*\in\mathbb{R}^d$ be a solution to
the population $K$-means problem, in the sense that they minimize
\[
R(\bolds{\mu}_1,\ldots,\bolds{\mu}_K)=\int\min
_{k=1,\ldots,K}\|\bx-\bolds\mu_k\|\,dP(\bx),
\]
where $\|\cdot\|$ denotes the usual Euclidean norm. Then, the $K$-means
algorithm assigns an arbitrary point
in $\mathbb R^d$ to the group whose center is closer, so that the ideal
population clustering is given
by $\mathscr C=\{C_1,\dots,C_K\}$, where
\[
C_k=\bigl\{\bx\in\mathbb R^d\colon\bigl\|\bx-\bolds{
\mu}_k^*\bigr\|\leq\bigl\|\bx -\bolds{\mu}_j^*\bigr\| \mbox{ for all }j
\neq k\bigr\}
\]
is the Voronoi cell corresponding to $\bolds\mu_k^*$, for
$k=1,\ldots,K$ (see \citep{GL00}, Chapter~4).
\begin{figure*}[t]

\includegraphics{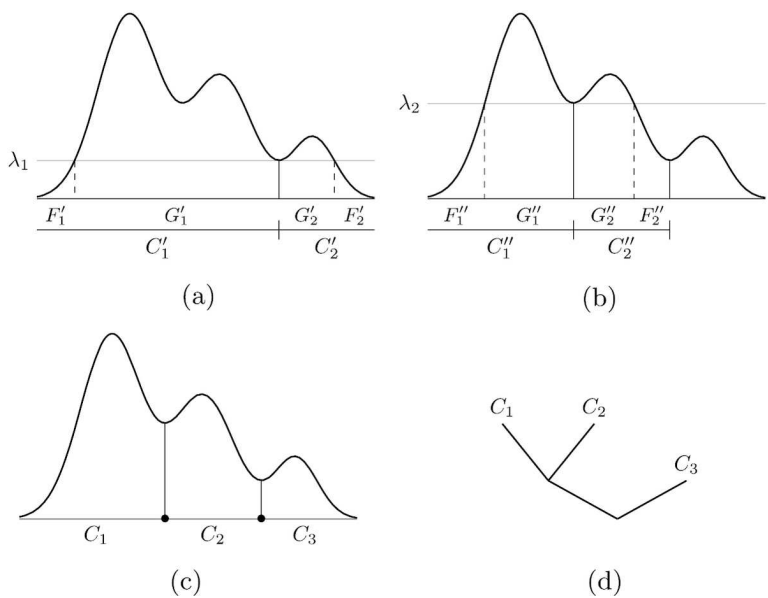}

\caption{Identification of clusters for the trimodal density example
using the cluster tree. Panel (a): first split; (b) second split; (c)
final partition; (d) cluster tree.} 
\label{fig2}
\end{figure*}

The ideal population clustering for mixture model clustering can be
derived in a similar way.
Assume that the underlying density is a mixture $f(\bx)=\sum_{k=1}^K\pi
_k  \cdot  f_k(\bx)$, where $\pi_k$
denotes the prior probability\vspace*{2pt} of the $k$th mixture component (with $\pi
_k>0$ and
$\sum_{k=1}^K\pi_k=1$), and $f_k(\bx)$ is the density of the $k$th
component. In this setup, assuming also that the mixture model is
identifiable, a point
$\bx\in\mathbb R^d$ is assigned to the group $k$ for which the a
posteriori probability $\pi_k
f_k(\bx)/f(\bx)$ is maximum, so the ideal population clustering that
$f$ induces has population
clusters
\[
C_k=\bigl\{\bx\in\mathbb R^d\colon\pi_k
f_k(\bx)\geq\pi_jf_j(\bx)\mbox{ for all }j
\neq k\bigr\}
\]
for $k=1,\ldots,K$. 



For the modal approach to clustering, however, the notion of ideal
population clustering is not so straightforward to formulate.
Informally, if the data-generating density $f$ has modes $\mathbf{M}_1,\ldots,\mathbf{M}_K$, then the population cluster $C_k$ is
defined as the domain of attraction of $\mathbf{M}_k$, for $k=1,\ldots
,K$. Most modal clustering algorithms are based on applying a
mode-seeking numerical method to the sample points and assigning the
same cluster to those data that are iteratively shifted to the same
limit value. Examples of such procedures include the mean shift
algorithm (\citep{FH75}), CLUES (\citep{WQZ07}) or the modal EM of
\citet{LRL07}, and further alternatives are described in a previous
unpublished version of this paper (\citep{C12}). Hence, from a practical
point of view, it is clear how a clustering of the data is constructed
on the basis of this notion of domain of attraction. The objective of
the next section is to describe in a precise way what is the population
goal that lies behind these algorithms. This aims to provide an answer,
in the case of modal clustering, to Question~1 in \citet{vLBD05}: ``How
does a desirable clustering look if we have complete knowledge about
our data generating process?''

\section{Describing the Population Goal of Modal Clustering Through
Morse Theory}\label{Morse}

The ideal population goal for modal clustering should reflect the
notion of a partition into regions of high density separated from each
other by regions of lower density. The following examples in one and
two dimensions are useful to illustrate the concept that we aim to formalize.

\begin{figure*}[b]

\includegraphics{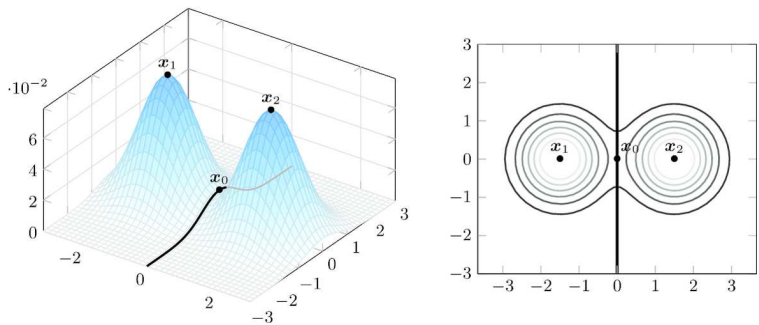}

\caption{Bidimensional example, with two groups clearly identifiable at
an intuitive level.}
\label{fig3}
\end{figure*}

In the one-dimensional case, it seems clear from Figure~\ref{fig2} how
this can be achieved. To
begin with, the level set methodology identifies the three clusters in
the density depicted in
Figure~\ref{fig1} by computing the cluster tree as described clearly
in \citet{NS10}: starting from
the $0$-level set, which corresponds to the whole real line in this
example (hence, it consists of a single connected component), $\lambda$
is increased until it reaches $\lambda_1$, where two components for the
$\lambda_1$-level set are found, $G_1'$ and $G_2'$, resulting in the cluster
tree splitting into two different branches [see Figure~\ref{fig2},
panel~(a)]. These two components
$G_1'$ and $G_2'$ are usually called cluster cores. They do not
constitute a clustering because there is some probability mass outside
$G_1'\cup G_2'$. But\vspace*{1pt} the
remaining parts $F_1'$ and $F_2'$, referred to as fluff in \citet{NS10},
can be assigned to either
the left or the right branch depending on whichever of them is closer.
Thus, at level $\lambda_1$
the partition $\mathbb R=C_1'\cup C_2'$ is obtained. The point dividing
the line into these two
components can be arbitrarily assigned to either of them; this
assignment makes no difference
because it leads to equivalent clusterings since a singleton has null
probability mass.

At level $\lambda_2$ the left branch $C_1'$ is further divided into two
branches [see panel (b) of Figure~\ref{fig2}]. Again, the two cluster
core components $G_1''$ and $G_2''$ do not form a partition of the set
$C_1'$ associated with the
previous node of the tree, but it is clear how the fluff $F_1''$ and
$F_2''$ can be assigned to form a
partition $C_1''\cup C_2''$ of $C_1'$. Since no further splitting of
the cluster tree is observed as
$\lambda$ increases, the final population clustering is $\{C_1'',
C_2'', C_2'\}$, renamed to $\{C_1,C_2,C_3\}$ in panel (c) of Figure~\ref{fig2}.

It is immediate to observe that the levels at which a connected
component breaks into two different
ones correspond precisely to local minima of the density function, so
an equivalent formulation
consists of defining population clusters as the connected components of
$\mathbb{R}$ minus the points where a local minimum is attained [the
solid circles in panel (c) of Figure~\ref{fig2}]. Notice that, unlike
the cluster tree, this definition does not involve the computation of
level sets for a range of levels, nor their cores and fluff, and in this
sense it constitutes a more straightforward approach to the very same
concept in the unidimensional setup.

\begin{figure*}[b]

\includegraphics{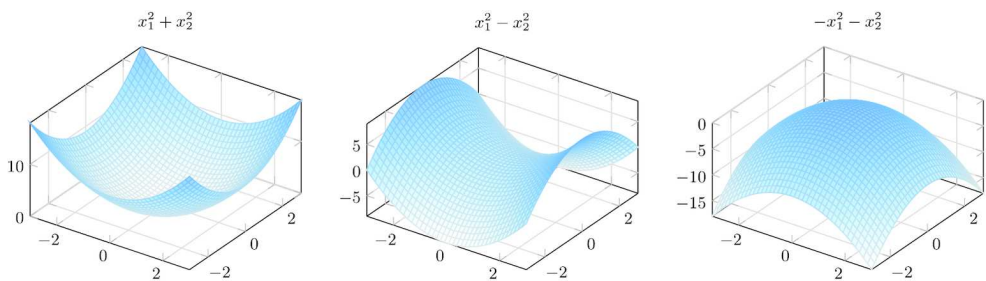}

\caption{The three possible configurations around a critical point of a
Morse function in the bidimensional case.}
\label{fig4}
\end{figure*}

To get an idea of how to generalize the previous approach to higher
dimensions, consider the
following extremely simple bidimensional example: an equal-proportion
mixture of two normal
distributions, each with identity variance matrix and centered at $\bmu_1=(-\frac{3}{2},0)$ and
$\bmu_2=-\bmu_1$, respectively. At an intuitive level, it is clear from
Figure~\ref{fig3} that the
most natural border to separate the two visible groups is the black
line. The problem is then: what is
exactly that line? Is it identifiable in terms of the features of the
density function in a
precise, unequivocal way? A nice way to answer these questions is by
means of Morse theory.

Morse theory is a branch of Differential Topology that provides tools
for analyzing the topology of
a manifold $M\subseteq\mathbb{R}^d$ by studying the critical points of a
smooth enough function
$f\colon M\to\mathbb{R}$. A classical reference book on this subject is
\citet{Mil63} and enjoyable
introductions to the topic can be found in \citet{Mat02} and
\citet{Jos11}, Chapter~7. A useful
application of Morse theory is for terrain analysis, as nicely
developed in \citet{Vit10}. In
terrain analysis, a mountain range can be regarded as the graph of a
function $f\colon M\to\mathbb{R}$, representing the terrain elevation, over a terrain $M\subseteq
\mathbb R^2$, just as in the left
graphic of Figure~\ref{fig3}. The goal of terrain analysis is to
provide a partition of $M$
through watersheds indicating the different regions, or catchment
basins, where water flows under
the effect of gravity.

The fundamentals of Morse theory can be extremely summarized as
follows. A smooth enough function $f\colon
M\to\mathbb{R}$ is called a \textit{Morse function} if all its critical
points are nondegenerate. Precisely, for our purposes, $f$ can be
considered smooth enough if it is three times continuously
differentiable. Here, the critical points of $f$ are understood as
those $\bx_0\in M$ for which the gradient $\D
f(\bx_0)$ is null, and nondegeneracy means that the determinant of the
Hessian matrix $\mathsf{H}
f(\bx_0)$ is not zero. For such points the \textit{Morse index} $m(\bx_0)$
is defined as the number of
negative eigenvalues of $\mathsf{H}f(\bx_0)$.

Morse functions can be expressed in a fairly simple form in a
neighborhood of a critical point
$\bx_0$, as the result known as Morse lemma shows that it is possible
to find local coordinates
$x_1,\dots,x_n$ such that $f$ can be written as $f(\bx_0)\pm x_1^2\pm
\cdots\pm x_d^2$ around $\bx_0$,
where the number of minus signs in the previous expression is precisely
$m(\bx_0)$. For example,
for $d=2$ the three possible configurations for a critical point are
shown in Figure~\ref{fig4},
corresponding to a local minimum, a saddle point and a local maximum
(from left to right), with
Morse indexes 0, 1 and 2, respectively.

The decomposition of $M$ suggested by Morse theory is made in terms of
the unstable and/or stable manifolds of the critical points of $f$ as
explained next. Consider the initial value problem defined by the minus
gradient vector of a smooth enough function $f$. For a given value of
$\bx\in M$ at time $t=0$, the integral curve $\bnu_\bx\colon\mathbb
R\to M$ of such an initial value problem is the one satisfying
%
\begin{equation}
\label{ngflow} \bnu_\bx'(t)=-\D f \bigl(
\bnu_\bx(t) \bigr),\quad \bnu_\bx(0)=\bx
\end{equation}
and the set of all these integral curves is usually referred to as the
negative gradient flow. Since the minus gradient vector defines the
direction of steepest descent of $f$, these curves (or, properly
speaking, their images through $f$) represent the trajectories of the
water flow subject to gravity.

With respect to the negative gradient flow, the \textit{unstable manifold}
of a critical point $\bx_0$ is defined as the set of points whose
integral curve starts at $\bx_0$, that is,
\[
W^u_-(\bx_0)= \Bigl\{\bx\in M\colon\lim
_{t\to-\infty}\bnu_\bx(t)=\bx _0 \Bigr\}.
\]
Analogously, the stable manifold of $\bx_0$ is the set of points whose
integral curve finishes at $\bx_0$, that is, $W^s_-(\bx_0)= \{\bx\in
M\colon\lim_{t\to+\infty}\bnu_\bx(t)=\bx_0 \}$. It was first noted
by \citet{Thom49} that the class formed by the unstable manifolds
corresponding to all the critical points of $f$ provides a partition of
$M$ (the same is true for the stable manifolds). Furthermore, the
unstable manifold $W^u_-(\bx_0)$ has dimension $m(\bx_0)$.

The main contribution of this section is the definition of the
population modal clusters of a density $f$ as the unstable
manifolds of the negative gradient flow corresponding to local maxima
of $f$. That is, if $\mathbf{M}_1,\ldots,\mathbf{M}_K$ denote the
modes of $f$, then the ideal population goal for modal clustering is
$\mathscr{C}=\{C_1,\ldots,C_K\}$, where $C_k=W^u_-(\mathbf{M}_k)$, for
$k=1,\ldots, K$. Or in a more prosaic
way, in terms of water flows, a modal cluster is just the region of the
terrain that would be flooded by
a fountain emanating from a peak of the mountain range.

Although this is an admittedly cumbersome definition, going back to
Figure~\ref{fig3}, it is clear
that it just describes the notion that we were looking for. The
critical point $\bx_0=(0,0)$ is a
saddle point, thus having Morse index 1, and the black line is
precisely its associated unstable
manifold, $W^u_-(\bx_0)=\{0\}\times\mathbb{R}$, which is a manifold of
dimension 1. The remaining
two critical points are local maxima, and their respective unstable
manifolds are
$W^u_-(\bx_1)=(-\infty,0)\times\mathbb{R}$ and $W^u_-(\bx_2)=(0,\infty
)\times\mathbb{R}$, manifolds
of dimension 2 so that we can partition $\mathbb{R}^2=W^u_-(\bx_0)\cup
W^u_-(\bx_1)\cup
W^u_-(\bx_2)$, showing $W^u_-(\bx_1)$ and $W^u_-(\bx_2)$ as two
population clusters separated by
the border $W^u_-(\bx_0)$, which is a null-probability set.

Notice that this definition also applies to the previous univariate
example in Figure~\ref{fig2}: the clusters $C_1$, $C_2$ and $C_3$ are
just the unstable manifolds of the three local maxima (they are
manifolds of dimension 1), and for the two local minima their unstable
manifolds have dimension 0, so they include only the respective points
of local minima.

\begin{figure*}[t]

\includegraphics{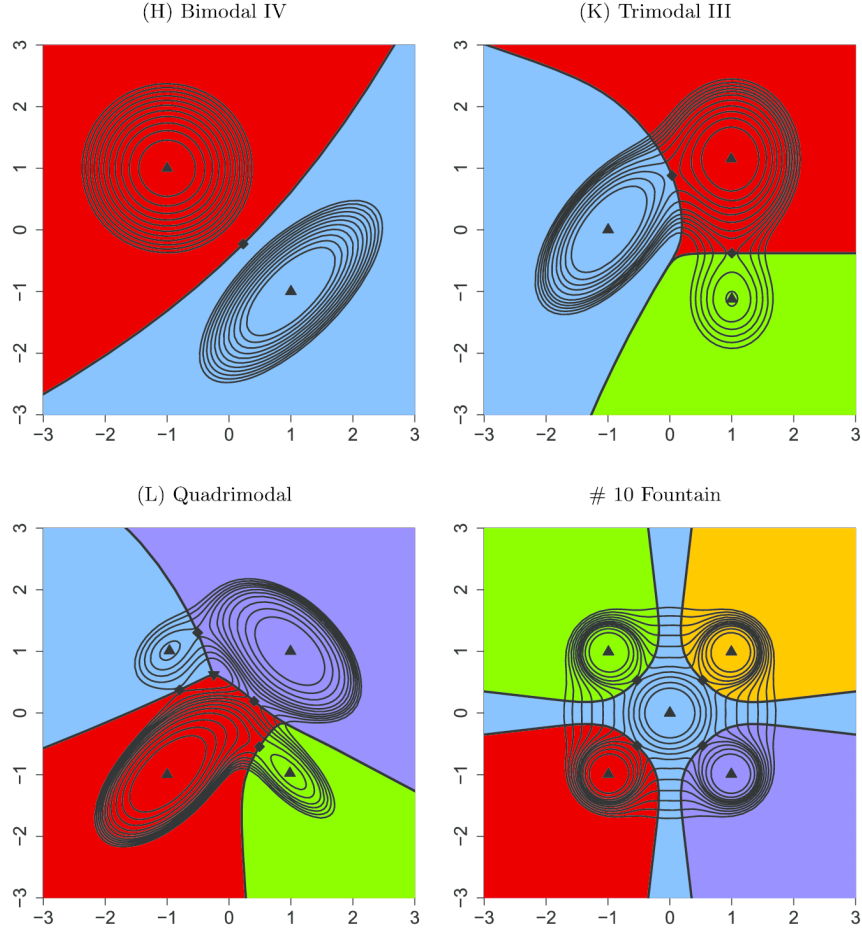}

\caption{Ideal modal population clustering for some normal mixtures densities.}
\label{examples}
\end{figure*}


Moreover, if we focus on the gradient flow, instead of the negative
gradient flow, then its integral curves satisfy
\begin{eqnarray*}
&& \bgamma_\bx'(t)=\D f \bigl(\bgamma_\bx(t)
\bigr),\quad \bgamma_\bx(0)=\bx;
\end{eqnarray*}
the unstable manifold for the negative gradient flow becomes the stable
manifold for the gradient
flow and viceversa. Therefore, we could equivalently define the cluster
associated to a mode
$\bx_0$ of the density as its stable manifold with respect to the
gradient flow, that is,
$W^s_+(\bx_0)=\{\bx\in M\colon\lim_{t\to\infty}\bgamma_\bx(t)=\bx_0
\}=W^u_-(\bx_0)$. This is a precise
formulation of the notion of domain of attraction of the mode $\bx_0$,
since $W^s_+(\bx_0)$
represents the set of all the points that climb to $\bx_0$ when they
follow the steepest ascent
path defined by the gradient direction. Moreover, estimating this path
is precisely the goal of the mean shift algorithm (see \citep{ACMP13}).

\subsection{Examples}

In Figure~\ref{examples} we give further examples of how the ideal
population goal of modal clustering looks for three of the bivariate
normal mixture densities included in \citet{WJ93}, namely, with their
terminology, densities (H) Bimodal IV, (K) Trimodal III and (L)
Quadrimodal, plus the normal mixture \#10 Fountain from \citet{Ch09}.
These densities have a number of modes ranging from two to five,
respectively, and hence that is the true number of population clusters
for each of these models, in the sense of modal clustering.

Each graph contains a contour plot of the density function; the
location of the modes is marked with a triangle pointing upward
($\blacktriangle$), the saddle points with a rotated square
($\blacklozenge$), and the only local minimum, appearing
in the plot of the Quadrimodal density, is marked with a triangle
pointing downward ($\blacktriangledown$). The thick lines passing
through the saddle points are their corresponding unstable manifolds
and represent the border between the different population clusters.

All these features have been computed numerically, making use of some
results from the thorough analysis of normal mixture densities given in
\citet{RL05}. For instance, the Newton--Raphson method has been used for
the location of the modes by finding a zero gradient point starting
from the component means, taking into account that both the location of
the modes and component means are different, but very close. Next, the
saddle points are searched along the ridgeline that connects every two
component means, since all the critical points of the density must lie
on this curve, by Theorem~1 in \citet{RL05}. Finally, the borders
between the population clusters are obtained by numerically solving the
initial value problem (\ref{ngflow}), starting from a point slightly
shifted from each saddle point, along the direction of the eigenvector
of its Hessian corresponding to a negative eigenvalue.

%
%

%

\section{Comparing Clusterings}\label{seccons}


Whatever the notion of ideal population clustering the researcher may
use, in practice, this population goal has to be approximated from the
data. Therefore, to evaluate the performance of a clustering method, it
is necessary to introduce a loss function to measure the distance
between a data-based clustering and the population goal or, more
generally, to have a notion of \mbox{distance} between two whole-space
clusterings. In this section, two proposals are derived by extending
two well-known notions of distance between sets to distances between
clusterings.


Recall that some clustering methods do not produce a partition of the
whole feature space, but only a clustering of the data. A good deal of
measures to evaluate the distance between two clusterings of the data
have been proposed in the literature. The work of \citet{M07} provides
both a comprehensive survey of the most
used existing measures as well as a deep technical study of their main
properties, and, for instance,
\citet{AB73} or \citet{D81} include further alternatives. But it should
be stressed that all these proposals concern only
partitions of a finite set. Here, on the contrary, our interest lies on
developing two new notions of distance between whole-space clusterings.

Let $\mathscr{C}$ and $\mathscr{D}$ be two clusterings of a probability
distribution $P$, and assume
for the moment that both have the same number of clusters, say,
$\mathscr{C}=\{C_1,\ldots,C_r\}$ and
$\mathscr{D}=\{D_1,\ldots,D_r\}$. The first step to introduce a distance
between $\mathscr{C}$ and
$\mathscr D$ is to consider a distance between sets. Surely the two
distances between sets most
used in practice are the Hausdorff distance and the distance in
measure; see \citet{CF10}. The
Hausdorff distance is specially useful when dealing with compact sets
(it defines a metric in the
space of all compact sets of a metric space), as it tries to capture
the notion of physical proximity
between two sets (\citep{RC03}). In contrast, given a measure~$\mu$, the
distance in $\mu$-measure between two sets $C$ and $D$
refers to $\mu(C\triangle D)$, that is, to the content of their
symmetric difference $C\triangle D=(C\cap D^c)\cup(C^c\cap D)$. It
defines a metric on the set of all measurable subsets of a measure
space, once two sets
differing in a null-measure set are identified to be the same.

\subsection{A Distance in Measure Between Clusterings}

Although we will return to the Hausdorff distance later, our first
approach to the notion
of distance between $\mathscr{C}$ and $\mathscr{D}$ relies primarily on
the concept of distance in $\mu$-measure, and the measure involved is
precisely the probability measure $P$. From a practical point of view,
it does not seem so important that the clusters of a data-based
partition get physically
close to those of the ideal clustering. Instead, it is desirable that
the points that are incorrectly assigned do not represent a very
significant portion of the distribution. This corresponds to the idea
of perceiving two clusters $C\in\mathscr{C}$ and $D\in\mathscr{D}$
(resulting from different clusterings) as close
when $P(C\triangle D)$ is low. In this sense, the closeness between $C$
and $D$ is quantified by their distance in $\mu$-measure for the
particular choice $\mu=P$.

Therefore, for two clusterings $\mathscr{C}$ and $\mathscr{D}$ with the
same number of clusters, a sensible notion of distance is obtained by adding
up the contributions of the pairwise distances between their components
once they have been
relabeled, so that every cluster in $\mathscr{C}$ is compared with its
most similar counterpart in
$\mathscr{D}$. In mathematical terms, the distance between $\mathscr{C}$
and $\mathscr{D}$ can be
measured by
%
\begin{equation}
\label{d1dist} d_1(\mathscr{C},\mathscr{D})=\min
_{\sigma\in\mathcal{P}_r}\sum_{i=1}^rP(C_i
\triangle D_{\sigma(i)}),
\end{equation}
where $\mathcal{P}_r$ denotes the set of permutations of $\{1,2,\break\ldots, r\}$.

It can be shown that $d_1$
defines a metric in the space of all the partitions with the same
number of components, once two
such partitions are identified to be the same if they differ only in a
relabeling of their
components. Moreover, the minimization problem in (\ref{d1dist}) is
usually known as the \textit{linear sum assignment problem} in the
literature of Combinatorial Optimization, and it represents a
particular case of the well-known Monge--Kantorovich transportation
problem. A comprehensive treatment of assignment problems can be found
in \citet{BDM09}.

If a partition is understood as a vector in the product space of
measurable sets, with the
components as its coordinates, then $d_1$ resembles the $L_1$ product
distance, only adapted to
take into account the possibility of relabeling the components. This
seems a logical choice given
the additive nature of measures, as it adds up the contribution of each
distance between the partition components as
described before. However, it would be equally possible to consider any
other $L_p$ distance,
leading to define
\[
d_p(\mathscr C,\mathscr D)=\min_{\sigma\in\mathcal P_r} \Biggl\{\sum
_{i=1}^rP(C_i\triangle
D_{\sigma(i)})^p \Biggr\}^{1/p}
\]
for $p\geq1$ and also $d_\infty(\mathscr{C},
\mathscr{D})=\min_{\sigma\in\mathcal{P}_r}\max\{P(C_i\cdot\triangle D_{\sigma
(i)})\colon i=1,\ldots,r\}$. The
minimization problem defining $d_\infty$ is also well known under the
name of the linear
bottleneck assignment problem, and its objective function is usually
employed if the interest is
to minimize the latest completion time in parallel computing
(see \citep{BDM09}, Section~6.2). Still, in
the context of clustering, surely the $d_1$ distance seems the most
natural choice among all the
$d_p$ possibilities, due to its clear interpretation.

\begin{figure}

\includegraphics{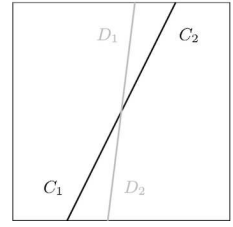}

\caption{When computing the distance $d_1$ between the two clusterings
$\mathscr{C}=\{C_1,C_2\}$ (black) and $\mathscr{D}=\{D_1,D_2\}$ (grey),
it is found that $C_1\cap D_1^c=C_2^c\cap D_2$ and $C_1^c\cap
D_1=C_2\cap D_2^c$, so the content of each of these two discrepancy
regions is added twice in $d_1(\mathscr{C},\mathscr{D})$.}
\label{fig5a}
\end{figure}

Nevertheless, the definition of the $d_1$ distance involves some kind
of redundancy, due to
the fact that $\mathscr{C}$ and $\mathscr{D}$ are (essential) partitions
of $\mathbb{R}^d$, because the two
disjoint sets that form every symmetric difference in fact appear twice
in each of the sums in
(\ref{d1dist}); see Figure~\ref{fig5a}.
%
More precisely, taking into account that $P(C\triangle
D)=P(C)+P(D)-2P(C\cap D)$, it follows that for every $\sigma\in\mathcal{P}_r$
%
\begin{eqnarray}
\quad\sum_{i=1}^rP(C_i
\triangle D_{\sigma(i)}) &=& 2-2\sum_{i=1}^r
P(C_i\cap D_{\sigma(i)})
\nonumber
\\[-8pt]
\label{fact2}
\\[-8pt]
\nonumber
&=& 2 P \Biggl( \Biggl\{\bigcup
_{i=1}^r(C_i\cap D_{\sigma
(i)})
\Biggr\}^c \Biggr).\hspace*{-10pt}
\end{eqnarray}
To avoid this redundancy, our eventual suggestion to measure the
distance between $\mathscr C$ and
$\mathscr D$, based on the set distance in $P$-measure, is
$d_P(\mathscr{C},\mathscr{D})=\frac{1}2 d_1(\mathscr{C},\mathscr{D})$.

If the partitions $\mathscr{C}$ and $\mathscr{D}$ do not have the same
number of clusters, then as many empty set components as needed are
added so that both partitions include the same number of components, as
in \citet{CDGH06}, and the distance between the extended partitions is
computed as before. Explicitly, if $\mathscr{C}=\{C_1,\ldots,C_r\}$ and
$\mathscr{D}=\{D_1,\ldots,D_s\}$ with $r<s$, then, writing
$C_i=\varnothing$ for $i=r+1,\ldots,s$, we set
\begin{eqnarray*}
&& d_P(\mathscr{C},\mathscr{D})\\[-2pt]
&&\quad=\frac{1}2\min
_{\sigma\in\mathcal{P}_s}\sum_{i=1}^sP(C_i
\triangle D_{\sigma(i)})\\[-2pt]
&&\quad=\frac{1}2\min_{\sigma\in\mathcal{P}_s} \Biggl\{
\sum_{i=1}^rP(C_i\triangle
D_{\sigma(i)})+\sum_{i=r+1}^sP(D_{\sigma(i)})
\Biggr\}.
\end{eqnarray*}
Thus, the term $\sum_{i=r+1}^sP(D_{\sigma(i)})$ can be interpreted as a
penalization for unmatched probability mass.\looseness=1
\begin{figure}

\includegraphics{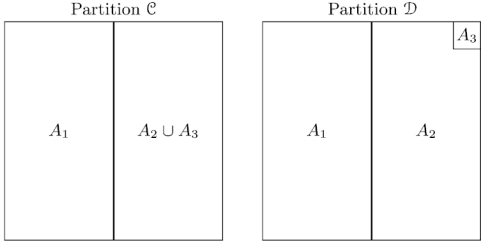}

\caption{Two partitions of the unit square that do not differ much if
$A_3$ has low probability.}
\label{fig5}
\end{figure}

The idea is that two partitions such as those shown in Figure~\ref{fig5} do not differ much if $A_3$
has low probability, even if they do not have the same number of
clusters. For the partitions in
Figure~\ref{fig5}, denote $\mathscr{C}=\{C_1,C_2\}$ and $\mathscr D=\{
D_1,D_2,D_3\}$ with
$C_1=D_1=A_1$, $C_2=A_2\cup A_3$, $D_2=A_2$, $D_3=A_3$, and assume that
$P(A_1)=0.5$, $P(A_2)=0.45$
and $P(A_3)=0.05$. Then, it can be shown that $d_P(\mathscr{C},
\mathscr{D})=0.05$. In (\ref{d1dist})
every cluster of $\mathscr{C}$ is matched to some cluster in
$\mathscr{D}$, depending on the
permutation for which the minimum is achieved. When
$\mathscr{C}$ has
less clusters than $\mathscr{D}$, some of the components of $\mathscr{D}$ will be matched with the
empty set, indicating that they
do not have an obvious match in $\mathscr{C}$s or that they are
unimportant. In the previous
example, the minimum is achieved when $C_1$ is matched with $D_1$,
$C_2$ with $D_2$ and $D_3$ is
matched with the empty set.

Indeed, if the existence of unmatched probability mass is considered to
be of greater concern, it is always possible to modify the distance in
$P$-measure by introducing a tuning parameter $\lambda\geq0$ to assign
a different weight to the penalization, thus mimicking other existing
procedures as penalized regression or pruning of decision trees. In
this case, the distance would be defined as
\begin{eqnarray*}
&& d_{P,\lambda}(\mathscr{C},\mathscr{D})\\[-2pt]
&&\quad=\frac{1}2\min
_{\sigma\in\mathcal{P}_s} \Biggl\{\sum_{i=1}^rP(C_i
\triangle D_{\sigma(i)})+\lambda\sum_{i=r+1}^sP(D_{\sigma(i)})
\Biggr\},
\end{eqnarray*}
so that $d_P(\mathscr{C},\mathscr{D})=d_{P,1}(\mathscr{C},\mathscr{D})$.

It is interesting to note that $d_P(\mathscr{C},\mathscr{D})$ can be
estimated in a natural way by
replacing $P$ with the empirical measure based on the data $\bX_1,\dots
,\bX_n$, leading to
\begin{eqnarray*}
&&\widehat{d_P}(\mathscr{C},\mathscr{D})\\
&&\quad=\frac{1}{2n}\min
_{\sigma\in\mathcal{P}_s} \Biggl\{\sum_{i=1}^r
\sum_{j=1}^nI_{C_i\triangle D_{\sigma(i)}}(\bX
_j)\\
&&\qquad{}+\sum_{i=r+1}^s\sum
_{j=1}^nI_{D_{\sigma(i)}}(\bX_j) \Biggr\},
\end{eqnarray*}
where $I_A$ denotes the indicator function of the set $A$. When $r=s$,
it follows from (\ref{fact2}) that an alternative expression for
$d_P(\mathscr{C},\mathscr{D})$ is
\[
d_P(\mathscr{C},\mathscr{D})=1-\max_{\sigma\in\mathcal{P}_r}\sum
_{i=1}^rP(C_i\cap
D_{\sigma(i)})
\]
and, therefore, its sample analogue,
\[
\widehat{d_P}(\mathscr{C},\mathscr{D})=1-\frac{1}n\max
_{\sigma\in\mathcal{P}_r}\sum_{i=1}^r\sum
_{j=1}^nI_{C_i\cap D_{\sigma(i)}}(
\bX_j),
\]
coincides with the so-called classification distance between two
clusterings of the data, whose
properties are explored in \citeauthor{M05}  (\citeyear{M05}, \citeyear{M07}, \citeyear{M12}). For $r<s$, however,
$\widehat{d_P}$ differs from
the classification distance (which does not include the penalty term),
but it corresponds exactly
with the transfer distance, studied in detail in
\citet{CDGH06} (see also \citep{D08}). Extending
the properties of the transfer distance to its population counterpart
suggests an interpretation of
$d_P(\mathscr{C},\mathscr{D})$ as the minimal probability mass that needs
to be moved to transform
the partition $\mathscr{C}$ into $\mathscr{D}$, hence the connection with
the optimal transportation
problem.

The above argument allows to recognize $d_P(\mathscr{C},\mathscr{D})$ as
the population version of some commonly used empirical distances
between partitions of a data set. However, it should be noted that the
estimate $\widehat{d_P}(\mathscr{C},\mathscr{D})$ requires the two
clusterings to be fully known and, hence, it may not be very useful if
the goal is to approximate the distance between the ideal population
clustering and a data-based clustering.

\subsection{A Hausdorff Distance Between Clusterings}

An alternative notion of distance between two clusterings based on the
Hausdorff metric has been kindly suggested by Professor Antonio Cuevas,
noting that precisely this distance was used in \citet{P81} to measure
the discrepancy between the set of sample $K$-means and the set of
population $K$-means. If $(X,\rho)$ is a metric space and $A,B\subseteq
X$ are two nonempty subsets of $X$, the Hausdorff distance between $A$
and $B$ is defined as
\[
d_H(A,B)=\max \Bigl\{\sup_{a\in A}\inf
_{b\in B}\rho(a,b),\sup_{b\in
B}\inf
_{a\in A}\rho(a,b) \Bigr\}
\]
or, equivalently, as
\[
d_H(A,B)=\inf\bigl\{\varepsilon>0\colon A\subseteq B^{\varepsilon}
\mbox{ and }B\subseteq A^{\varepsilon} \bigr\},
\]
where $A^\varepsilon=\bigcup_{a\in A}\{x\in X\colon\rho(x,a)\leq
\varepsilon\}$, and $B^\varepsilon$ is defined analogously.

In the context of clustering, $X$ can be taken to be the metric space
consisting of all the sets of $\mathbb{R}^d$ equipped with the distance
$\rho(C,D)=P(C\triangle D)$, once two sets with $P$-null symmetric
difference have been identified to be the same. Then any two
clusterings $\mathscr{C}=\{C_1,\ldots,C_r\}$ and $\mathscr{D}=\{D_1,\ldots
,D_s\}$ can be viewed as (finite) subsets of $X$ and, therefore, the
Hausdorff distance between $\mathscr{C}$ and $\mathscr{D}$ is defined as
\begin{eqnarray*}
&& d_H(\mathscr C,\mathscr D)\\
&&\quad= \max \Bigl\{\max_{i=1,\ldots,r}
\min_{j=1,\ldots
,s}P(C_i\triangle D_j),\\
&&\qquad\max
_{j=1,\ldots,s}\min_{i=1,\ldots
,r}P(C_i\triangle
D_j) \Bigr\}
\\
&&\quad= \inf\bigl\{\varepsilon>0\colon\mathscr C\subseteq\mathscr D^{\varepsilon
}
\mbox{ and }\mathscr D\subseteq\mathscr C^{\varepsilon} \bigr\}.
\end{eqnarray*}
To express it in words, $d_H(\mathscr{C},\mathscr{D})\leq\varepsilon$
whenever for every $C_i\in\mathscr{C}$ there is some $D_j\in\mathscr{D}$
such that $P(C_i\cdot\triangle D_j)\leq\varepsilon$ and vice versa. Hence,
as noted by \citet{P81}, if $\varepsilon$ is taken to be less than one
half of the minimum of distance between the clusters within $\mathscr{C}$ and also less than one half of the minimum distance between the
clusters within $\mathscr{D}$, then $d_H(\mathscr{C},\mathscr{D})\leq
\varepsilon$ implies that $\mathscr{C}$ and $\mathscr{D}$ must
necessarily have the same number of clusters.

The Hausdorff distance can be regarded as a uniform distance between
sets. It is not hard to show, using standard techniques from the Theory
of Normed Spaces, that when $r=s$ we have
\[
d_H(\mathscr{C},\mathscr{D})\leq2 d_P(\mathscr{C},
\mathscr{D})\leq r d_H(\mathscr{C},\mathscr{D}).
\]
However, when $r<s$ the distance $d_H$ can be more demanding than
$d_P$, meaning that both
partitions have to be really close so that their Hausdorff distance
results in a small value. For
instance, it can be checked that for the two clusterings of the
previous example, shown in Figure~\ref{fig5}, the Hausdorff distance between them is $d_H(\mathscr{C},\mathscr{D})=0.45$, mainly due
to the fact that $C_2$ and $D_3$ are far from each other, since
$P(C_2\triangle D_3)=P(A_2)=0.45$.

A clear picture of the difference between $d_H$ and $d_P$ is obtained
by arranging all the
component-wise distances $P(C_i\triangle D_j)$ into an $r\times s$
matrix. Then, the \mbox{Hausdorff}
distance is obtained by computing all the row-wise and column-wise
minima and taking the maximum of
all of them. In contrast, for the distance in $P$-measure the first
step when $r<s$ is to add $s-r$
row copies of the vector $(P(D_1),\ldots,P(D_s))$ to the matrix of
component-wise distances, and
then compute the distance in $P$-measure as half the minimum possible
sum obtained by adding up a
different element in each row. As a further difference, note that the
Hausdorff distance does not
involve a matching problem; instead, this distance is solely determined
by the two components that
are furthest from each other.

Obviously, a sample analogue is also obtained in this case by replacing
$P$ for the empirical
probability measure, leading to
\begin{eqnarray*}
&& \widehat{d}_H(\mathscr{C},\mathscr{D})\\
&&\quad=\frac{1}n\max
\Biggl\{\max_{i=1,\ldots
,r}\min_{j=1,\ldots,s}\sum
_{k=1}^nI_{C_i\triangle D_j}(\bX_k),\\
&&\qquad \max
_{j=1,\ldots,s}\min_{i=1,\ldots,r}\sum
_{k=1}^nI_{C_i\triangle D_j}(\bX _k) \Biggr
\},
\end{eqnarray*}
which seems not to have been considered previously as a distance
between two clusterings of the
data.

\subsection{Consistency of Data-Based Clusterings}

As indicated above, a data-based clustering is understood as any
procedure that induces a clustering
$\widehat{\mathscr{C}}_n$ of a probability distribution $P$ based on the
information obtained from a sample
$\bX_1,\ldots,\bX_n$ from $P$. Once a clustering methodology has been
chosen, and its ideal population goal\vspace*{1.5pt} $\mathscr{C}_0$ is clearly
identified, a data-based clustering $\widehat{\mathscr{C}}_n$ can be
said to be consistent if it gets closer to $\mathscr{C}_0$ as the sample
size increases. Formally, if
$d(\widehat{\mathscr{C}}_n,\mathscr{C}_0)\to0$ as $n\to\infty$ for some
of the modes of stochastic
convergence (in probability, almost surely, etc.), $d$ represents one of
the distances between
clusterings defined above or any other sensible alternative. Note that
a different notion of consistency, specifically intended for the
cluster tree approach, is studied in \citet{CD10}.

For density-based clustering, a plug-in strategy to obtain data-based
clusterings would consist of replacing the unknown density $f$ with
an estimator $\hat{f}_n$. Obvious candidates for the role of $\hat{f}_n$
include nonparametric density estimators for modal clustering or
\mbox{mixture} model density estimators with parameters fitted by maximum
likelihood for mixture model clustering. This is a very simple approach
that involves to some extent estimating the density function to solve
the clustering problem (unsupervised learning).

According to \citet{vL04}, page 21, this plug-in strategy may
not be a good idea because density estimation is a very difficult
problem, especially in high
dimensions. However, a similar situation is found in the study of
classification (supervised learning),
where the optimal classifier, the Bayes rule, depends on the regression
function of the random
labels over the covariates. Here, even if classification can be proved
to be a problem easier than
regression, nevertheless, regression-based algorithms for
classification play an important role in
the development of supervised learning theory (see  \citep{DGL96}, Chapter~6).

\begin{figure}[t]

\includegraphics{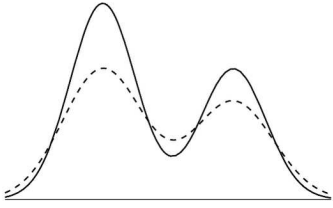}

\caption{Two density functions that are not close but induce exactly
the same clustering.}
\label{fig6}
\end{figure}

Along the same lines, Figure~\ref{fig6} illustrates why we should not
completely discard density estimation as an intermediate step for
clustering. Figure~\ref{fig6} shows a typical situation where the
solid line is the true density and the dashed line is a kernel density
estimator, since an expansion of its pointwise bias shows that, on
average, the kernel estimator underestimates the maxima and
overestimates the minima (\citep{WJ95}, page 21). But even if the two
density functions are not really close in any global sense, they
produce exactly the same clusterings of $\mathbb{R}$.

In any case, the following result shows that the plug-in strategy leads
to consistent data-based modal clusterings as long as the first and
second derivatives of the sequence of density estimators converge
uniformly to their true density counterparts.

\begin{theorem}\label{thmcons}
Let a Morse function $f$ be the density of a univariate probability
distribution $P$ with compact support, and denote by $\mathscr{C}_0$ the
ideal modal clustering that it induces, as defined in Section~\ref{Morse}. Let $\{\hat f_n\}$ be a sequence of density estimators\vspace*{1pt} such
that\vspace*{1pt} $\hat f^{(j)}_n\to f^{(j)}$ uniformly almost surely\vspace*{1pt} for $j=1,2$,
with $^{(j)}$ standing for the $j$th derivative. Denote by $\widehat
{\mathscr{C}}_n$ the modal clustering induced by $\hat f_n$. Then:
\begin{longlist}[(a)]
\item[(a)] $\#\widehat{\mathscr{C}}_n\to\#{\mathscr{C}}_0$ with
probability one as $n\to\infty$, where $\# A$ denotes the number of
elements in a set $A$.

%
\item[(b)] Both $d_P(\widehat{\mathscr{C}}_n,{\mathscr{C}}_0)\to0$ and
$d_H(\widehat{\mathscr{C}}_n,{\mathscr{C}}_0)\to0$ with probability one
as $n\to\infty$.
\end{longlist}
\end{theorem}

The proof of this result is shown in the
\hyperref[app]{Appendix}. The analysis of the
proposed distances between clusterings is greatly simplified in the
univariate case since the cluster boundaries are solely determined by
the points of local minima of the density. The extension of this result
for dimension $d\geq2$ seems quite a challenging open problem, since
the cluster boundaries in dimension $d$ are $(d-1)$-dimensional
manifolds which may have very intricate forms.

Part (a) shows that the number of clusters in $\widehat{\mathscr{C}}_n$ converges to the true number of clusters in ${\mathscr{C}}_0$
almost surely. As indicated in \citet{CFF00},
since $\#\widehat{\mathscr{C}}_n$ and $\#{\mathscr{C}}_0$ are integer-valued, this convergence is
equivalent to the fact that the event
\begin{eqnarray*}
&& \{\mbox{There exists } n_0\in\mathbb N\mbox{ such that }
 \\
&&\quad\#\widehat {\mathscr{C}}_n=\#{\mathscr{C}}_0\mbox{ for all } n\geq
n_0\}
\end{eqnarray*}
has probability one.

Note also that if $f^{(2)}$ is uniformly continuous and $\hat f_n$ are
kernel estimators with bandwidth $h=h_n$ based on a sufficiently
regular kernel, \citet{Sil78}, Theorem C, showed that a necessary and
sufficient condition for the uniform convergence condition in the
previous theorem to hold is just that $h\to0$ and $nh^5/\log n\to\infty
$ as $n\to\infty$ (see also \citep{Deh74} and \citep{BR78}).

\subsection{Asymptotic Loss Approximations}

The proof of Theorem~\ref{thmcons} reveals that, for big enough $n$,
the distance in measure and the Hausdorff distance between $\widehat
{\mathscr{C}}_n$ and $\mathscr C_0$ can be written as
\begin{eqnarray*}
d_P (\widehat{\mathscr{C}}_n,\mathscr{C}_0
)&=&\sum_{j=1}^{r-1} \bigl|F(\hat
m_{n,j})-F(m_j) \bigr| \quad\mbox{and}
\\
d_H (\widehat{\mathscr{C}}_n,\mathscr C_0
)&=&\max_{j=1,\ldots
,r-1} \bigl|F(\hat m_{n,j})-F(m_j) \bigr|,
\end{eqnarray*}
where $F$ is the distribution function of $P$. Here, $m_1,\dots
,m_{r-1}$ and $\hat m_{n,1},\dots,\hat m_{n,r-1}$ denote the local
minima of $f$ and $\hat{f}_n$, respectively\vspace*{1.5pt} (i.e., the cluster
boundaries of $\mathscr{C}_0$ and $\widehat{\mathscr{C}}_n$). From these
expressions the $L_1$ and $L_\infty$ nature of $d_P$ and $d_H$ is even
more clear.

Furthermore, under the conditions of Theorem~\ref{thmcons}, after two
Taylor expansions it is possible to obtain the approximations
\begin{eqnarray*}
\bigl|F(\hat m_{n,j})-F(m_j) \bigl| &\simeq &  f(m_j) \bigl|\hat
m_{n,j}-m_j \bigr|\\
& \simeq & \frac{f(m_j)}{f''(m_j)} \bigl|\hat
f_n'(m_j) \bigr|.
\end{eqnarray*}
This shows how not only the performance of $\widehat{\mathscr{C}}_n$ is
closely connected to the problem of first-derivative estimation, but
also that modal clustering is more difficult, as the density at the
cluster boundaries is higher and/or flatter as the intuition dictates.

In the case of kernel estimators, Proposition~4.1 of \citet{Ro88}
provides a precise description of the asymptotic behavior of $\hat
f_n'(m_j)$. Precisely, under some smoothness conditions it can be shown
that assuming that the bandwidth further satisfies $nh^7\to\beta^2$
with $0\leq\beta<\infty$, then $\hat f_n'(m_j)$ admits the representation
\[
\hat{f}_n'(m_j)=\bigl(nh^3
\bigr)^{-1/2}\sigma Z_n+\beta\mu
\]
for some explicit constants $\sigma>0$ and $\mu\in\mathbb R$, where
$Z_n$ is a sequence of asymptotically $N(0,1)$ random variables. This
representation could be helpful as a starting point to tackle the
problem of optimal bandwidth choice for kernel clustering, which has
only been treated briefly in the previous literature
(e.g., \citep{E11}, \citep{CD13}, \citep{ChM13}) and surely deserves further investigation.
However, we will not pursue this further here.

\section{Discussion}

At the time of comparing different clustering procedures, it is
necessary to have a ``ground truth,''
or population goal, that represents the ideal clustering to which the
clustering algorithms should
try to get close. The importance of having a clear population goal for
clustering is nicely highlighted in \citet{K09}, Chapter~8. Sometimes
this ideal population clustering is not so easy to specify, and of
course it
depends on the notion of cluster in which the researcher is interested.

Whereas the population goal is clearly defined for some clustering
methods, like $K$-means clustering or mixture model clustering, it
remained less obvious for modal clustering. Here, the ideal population
goal of modal clustering is accurately identified, making use of some
tools from Morse theory as the partition of the space induced by the
domains of attraction of the local maxima of the density function.

This definition of the modal clusters needs the probability density to
be smooth to a certain degree, specifically it
must be a 3-times continuously differentiable Morse function. It would
be appealing to extend this notion to density functions that are
not Morse functions, meaning either that they are smooth but have
degenerate critical points or
even that they are not differentiable to such extent. To treat the
first case, it might be useful to resort to the
theory of singularities of differential mappings, which is exhaustively
covered in the book by
\citet{AGLV98}, for instance. On the other hand, the study of the
nonsmooth case might start from
\citet{APS97}, where Morse theory for piecewise smooth functions is
presented. Here, the key role
would be played by the subgradient, which generalizes the concept of
the gradient for nonsmooth
functions.

Alternatively, as in \citet{D88}, a nonsmooth density $f$ could be
convolved with a mollifier $\phi_h$ to obtain a smoother version $\phi
_h*f$, so that the population modal clustering $\mathscr C_h$ of $\phi
_h*f$ is determined as in the smooth case, and then define the
population modal clustering of $f$ as the limit (in some sense) of
$\mathscr C_h$ as $h\to0$. Of course, further investigation on how to
properly formalize this notion would be required.


Once a clustering methodology with a clearly defined population goal
has been chosen, it is necessary to have a distance to measure the
accuracy of data-based clusterings as approximations of the ideal goal.
A second contribution of this paper is the introduction of two new loss
functions for this aim, which are valid for any clustering methodology.
Particularly, when applied to modal clustering, it is shown that the
plug-in approach leads to clustering consistency under mild assumptions.

A further interesting challenge for future research consists of
studying the choice of the parameters for the density estimators (the
bandwidth for kernel estimators, the mixture parameters for mixture
model estimators) that minimize the distance between the corresponding
data-based clustering and the true population clustering, as measured
by any of the distances between clusterings discussed in Section~\ref{seccons}. Or, maybe even better, to develop methods aimed to perform
modal clustering that do not necessarily rely on a pilot density
estimate, perhaps by somehow adapting those classification methods
whose construction is not based on a regression estimate.

\begin{appendix}\label{app}

\section*{Appendix: Proof of the Consistency~Theorem}

The proof uses some arguments from Theorem~3 in \citet{CGM91}; see also
Lemma~3 in \citet{GPVW13}.

First, since $f$ is a Morse function with compact support, it has only
finitely many isolated critical points (\citep{Mat02}, Corollary~2.19).
Assume that $f$ has $r$ local maxima and let $m_1<\cdots<m_{r-1}$
denote the local minima of $f$ so that the modal population clustering
induced by $f$ is defined as $\mathscr C_0=\{C_1,\ldots,C_r\}$ with
$C_j=(m_{j-1},m_j)$ for $j=1,\ldots,r$, where $m_0=-\infty$ and
$m_r=\infty$ (if $f$ has no local minimum, then $r=1$ and $\mathscr{C}_0=\{C_1\}=\{\mathbb{R}\}$).

We claim\vspace*{1pt} the following: with probability one, there exists $n_0\in
\mathbb N$ such that
$\hat f_n$ has exactly $r-1$ local minima for all $n\geq n_0$;
moreover, there exists $\varepsilon>0$ such that every $\hat f_n$ with
$n\geq n_0$ has exactly one local minimum $\hat m_{n,j}$ in
$[m_j-\varepsilon,m_j+\varepsilon]$ for all $j=1,\ldots,r-1$. To prove
this claim, notice that since $f''(m_j)>0$ for all $j$, and $f''$ is
continuous, it is possible to find some $\varepsilon>0$ such that
$f''(x)>0$ on $[m_j-\varepsilon,m_j+\varepsilon]$, for all~$j$. The
almost sure uniform convergence of $\hat f_n''$ to~$f''$ implies that
there is some $n_0\in\mathbb{N}$ such that, with a possibly smaller
$\varepsilon$, all $\hat f_n''$ with $n\geq n_0$ are strictly positive
on those intervals as well. On the other hand, on each of these
intervals $f'$ is strictly increasing and since $f'(m_j)=0$, it must go
from negative to positive. But the uniform convergence of $\hat f_n'$
to $f'$ implies that also $\hat f_n'$ must go from negative to positive
(perhaps with a smaller $\varepsilon$) for big enough $n$. Therefore,
all of them must have a critical point there,\vspace*{1pt} and since we previously
showed that $\hat{f}_n''>0$, this means both that the critical point is
a local minimum and that there cannot be any more of them in such
neighborhoods of the local minima. A similar argument shows that, for
big enough $n$, all the $\hat f_n$ with $n\geq n_0$ must also have a
local maximum in a small enough neighborhood around the modes of $f$,
and that there cannot be other critical points of $\hat{f}_n$ outside
these neighborhoods.

Furthermore, using standard arguments in $M$-estimation theory, under
these conditions it follows that also $\hat m_{n,j}$ converges to $m_j$
as $n\to\infty$: to show this, notice that given an arbitrary $\eta>0$,
small enough so that $\eta<\varepsilon$, the value of $\delta:=\inf\{
|f'(x)|\colon\eta\leq|x-m_j|\leq\varepsilon\}$ is strictly positive.\vspace*{-1pt}
Hence, from the almost sure uniform convergence $\hat f_n'\to f'$ it
follows that, with probability one, for all big enough $n$ we have
$|\hat f_n'(x)|>\delta/2>0$ whenever $\eta\leq|x-m_j|\leq\varepsilon$.
Since $|\hat m_{n,j}-m_j|\leq\varepsilon$ and $\hat f_n'(\hat
m_{n,j})=0$, this implies that $|\hat m_{n,j}-m_j|<\eta$.

In this situation, for the clustering $\widehat{\mathscr C}_n=\{
\widehat C_{n1},\ldots,\break\widehat C_{n,r}\}$ induced by $\hat f_n$ [with
$\widehat C_{n,j}=(\hat m_{n,j-1},\hat m_{n,j})$, $\hat m_{n,0}=-\infty
$ and $\hat
m_{n,r}=\infty$], taking a small enough $\varepsilon$, the distance in
$P$-measure and the Hausdorff distance between $\widehat{\mathscr C}_n$
and $\mathscr C$ can be simply written as
\begin{eqnarray*}
d_P (\widehat{\mathscr{C}}_n,\mathscr{C}_0
)&=&\sum_{j=1}^{r-1} \bigl|F(\hat
m_{n,j})-F(m_j) \bigr|,
\\
d_H (\widehat{\mathscr{C}}_n,\mathscr{C}_0
)&=&\max_{j=1,\ldots
,r-1} \bigl|F(\hat m_{n,j})-F(m_j) \bigr|,
\end{eqnarray*}
respectively, where $F$ is the distribution function of $P$. Therefore,
the convergence of the estimated local minima to the true local minima
of $f$ yields the result.
\end{appendix}

\section*{Acknowledgments}
The author wishes to thank Professor Antonio Cuevas from Universidad
Aut\'onoma de Madrid as well as Professor Ricardo Faro from Universidad
de Extremadura for insightful
conversations and suggestions concerning the material of Section~\ref{seccons}. The paper by
\citet{RL05}, in which interesting connections between Morse theory and
the topography of
multivariate normal mixtures are illustrated, was thought-provoking
enough to inspire part of this
paper.

Supported in part by Spanish Ministerio de Ciencia y
Tecnolog\'{\i}a
projects MTM2010-16660, MTM2010-17366 and MTM2013-44045-P, and by the
Gobierno de Extremadura Grant GR10064.





\end{document}